\documentclass[12pt]{amsart}

\usepackage{fullpage}
\usepackage{amssymb}
\usepackage{amsmath}
\usepackage{url}
\usepackage{color}
\usepackage[colorlinks, breaklinks, linkcolor=blue]{hyperref} 
\usepackage{breakurl}

\newcommand{\height}{\operatorname{ht}}
\newcommand{\maxim}{\operatorname{max}}

\newcommand\nn{\mathrm{n}}
\newcommand{\F}{\mathbb F}
\newcommand{\N}{\mathbb N}
\newcommand{\Z}{\mathbb Z}
\newcommand{\C}{\mathbb C}
\newcommand{\Q}{\mathbb Q}

\renewcommand{\c}{\mathfrak c}
\newcommand{\g}{\mathfrak g}
\newcommand{\m}{\mathfrak m}
\renewcommand{\u}{\mathfrak u}

\DeclareMathOperator{\Top}{top} 
\DeclareMathOperator{\GL}{GL} 
\DeclareMathOperator{\Mat}{Mat} 
 
\DeclareMathOperator{\U}{U}

\DeclareMathOperator{\modulo}{mod}

\DeclareMathOperator{\R}{R} 
\DeclareMathOperator{\I}{I} 
\DeclareMathOperator{\lcm}{lcm}

\numberwithin{equation}{section}

\newtheorem{thm}{Theorem}[section]
\newtheorem{lem}[thm]{Lemma}

\theoremstyle{definition}

\newtheorem{exmp}[thm]{Example}

\theoremstyle{remark}

\title[Calculating conjugacy classes in Sylow $p$-subgroups]
{Calculating conjugacy classes in Sylow $p$-subgroups \\ of finite
Chevalley groups of rank six and seven}


\author[S.~M.~Goodwin]{Simon M. Goodwin}
\address{School of Mathematics, University of Birmingham, Birmingham, B15
2TT, United Kingdom}  \email{s.m.goodwin@bham.ac.uk}

\author[P.~Mosch]{Peter Mosch} \author[G.\ R\"ohrle]{Gerhard R\"ohrle}
\address{Fakult\"at f\"ur Mathematik, Ruhr-Universit\"at Bochum, D-44780 Bochum, Germany}  \email{peter.mosch@rub.de} \email{gerhard.roehrle@rub.de}

 


\begin{document}

\begin{abstract}

Let $G(q)$ be a finite Chevalley group, where $q$ is a power of a good prime $p$, and let $U(q)$ be a Sylow $p$-subgroup of $G(q)$. Then a generalized version of a conjecture of Higman asserts that the number $k(U(q))$ of conjugacy classes in $U(q)$ is given by a polynomial in $q$ with integer coefficients. In \cite{GOROconj}, the first and the third authors developed an algorithm to calculate the values of $k(U(q))$. By implementing it into a computer program using {\sf GAP}, they were able to calculate $k(U(q))$ for $G$ of rank at most 5, thereby proving that for these cases $k(U(q))$ is given by a polynomial in $q$. In this paper we present some refinements and improvements of the algorithm that allow us to calculate the values of $k(U(q))$ for finite Chevalley groups of rank six and seven, except $E_7$. We observe that $k(U(q))$ is a polynomial, so that the generalized Higman conjecture holds for these groups. Moreover, if we write $k(U(q))$ as a polynomial in $q-1$, then the coefficients are non-negative. 

Under the assumption that $k(U(q))$ is a polynomial in $q-1$, we also give an explicit formula for the coefficients of $k(U(q))$ of degrees zero, one and two.

\end{abstract}

\maketitle

\section{Introduction}
\label{sec:intro}

Let $\GL_n(q)$ be the group of invertible $n\times n$ matrices with coefficients in the finite field $\F_q$, where $q$ is a power of the prime $p$, and let $\U_n(q)$ be the subgroup of unipotent upper triangular matrices. A well-known conjecture attributed to Higman is that the number $k(\U_n(q))$ of $\U_n(q)$-conjugacy classes in $\U_n(q)$ is given by a polynomial in $q$ with integer coefficients independent of $q$, see \cite{HI}. This conjecture has attracted the interest of many mathematicians including Thompson \cite{TH} and Robinson \cite{RB}. 

Using computer calculations, Vera-L\'opez and Arregi verified in \cite{VLAR} that the conjecture holds for $n\leq13$. The resulting polynomials have the additional property that, considered as polynomials in $q-1$, their coefficients are non-negative integers. We also note that Evseev has calculated these polynomials via an alternative approach, see \cite{EV}.  

Since the conjugacy classes of a finite group are in bijective correspondence with its complex irreducible characters, one can also approach the conjecture via character theory. This has been considered among others by Andr\'e \cite{AN}, Isaacs \cite{IScharupp} and Lehrer \cite{LH}. 

In this paper, we consider a generalization of Higman's conjecture. Let $G(q)$ be a finite Chevalley group, i.e.\ the group of $\F_q$-rational points of a simple algebraic group $G$ which is defined and split over $\F_q$. Assume that $p$ is good for $G(q)$ and let $U(q)$ be a Sylow $p$-subgroup of $G(q)$. Then the generalized conjecture claims that the number $k(U(q))$ of $U(q)$-conjugacy classes is given by a polynomial in $q$, and as a polynomial in $q-1$ it has non-negative integer coefficients. 

There has been a lot of interest recently in the conjugacy classes and the complex characters of $U(q)$, some of which gives evidence for (the generalized) Higman's conjecture.  For example, in \cite{ALP}, Alperin proved that the number $k(\U_n(q),\GL_n(q))$ of $\U_n(q)$-conjugacy classes in $\GL_n(q)$ is given by a polynomial in $q$ with integer coefficients. This was generalized by the first and the third authors in \cite{GOROflag} where they showed that $k(U(q),G(q))$ is a polynomial in $q$ when the centre of $G$ is connected and $G$ is not of type $E_8$. In case $G$ is of type $E_8$, the number $k(U(q),G(q))$ is given by one of two polynomials, depending on $q$ $\modulo$ 3.  It is conceivable that similar PORC (\textbf{P}olynomial \textbf{O}n \textbf{R}esidue \textbf{C}lasses) behaviour occurs for $k(U(q))$ if $G$ is of type $E_8$. For other recent developments, see for example \cite{GOROcount}, \cite{HH}, \cite{HLM} and \cite{LE}.

An algorithm was introduced in \cite{GOROconj} to calculate a parameterization of the
conjugacy classes of $U(q)$, and thus determine $k(U(q))$. In this paper we describe an improved version of this algorithm. Both versions of the program have been implemented in {\sf GAP} \cite{GAP}. The main idea of the algorithm (which is based on the results in \cite{GOconj}) is to replace the task of counting conjugacy classes by the geometric task of counting $\F_q$-rational points of quasi-affine varieties over finite fields, which parametrize the conjugacy classes of $U(q)$. The goal of the program is to determine these varieties and then calculate the number of $\F_q$-rational points from the polynomial equations which define them. In the previous version of the algorithm, it was necessary to inspect some output of the program and complete some calculations by hand and it was only possible to calculate $k(U(q))$ when the rank of $G$ is at most $5$. 

The improved version of the algorithm determines the polynomial equations much more effectively, and is significantly better at calculating the number of rational points in the varieties directly. With the aid of the improved version, we are able to prove the following theorem.

\begin{thm}
\label{thm:main}
Let $G$ be a split simple algebraic group defined over $\F_q$ of rank at most 7, excluding $E_7$, where $q$ is a power of a good prime $p$. Let $U$ be a maximal unipotent subgroup of $G$ which is also defined over $\F_q$. Then the number $k(U(q))$ of $U(q)$-conjugacy classes in $U(q)$ is given by a polynomial in $q$ with integer coefficients. Furthermore, if one considers $k(U(q))$ as a polynomial in $q-1$, then the coefficients are non-negative. 
\end{thm}

The theoretical background underlying our algorithm only holds for good primes. Therefore, the values for $k(U(q))$ given in this paper are only valid in this case. In \cite{BRGO}, Bradley and the first author calculated $k(U(q))$ for $q$ a power of a bad prime and $G$ of rank at most 4, excluding $F_4$. In these cases $k(U(q))$ is again given by polynomials with integer coefficients which are expectedly not the same as those for good $p$.  

As mentioned in \cite{GOROconj}, it is straightforward to adapt the program to calculate the number of $M(q)$-classes in $N_1(q)/N_2(q)$, where $M,N_1,N_2$ are normal unipotent subgroups of a Borel subgroup of $G$ defined over $\F_q$. For example, it is possible to calculate the number $k(U(q),U^{(l)}(q))$ of $U(q)$-conjugacy classes in the $l$-th term of the descending central series of $U(q)$ for $l\in\N$. We have done this for some cases where $G$ is of type $E_7$ and $E_8$. Here we also found that all values which we calculated were given by polynomials in $q$ with integer coefficients. While the previous version of the algorithm could calculate $k(U(q),U^{(l)}(q))$ for $G$ of type $E_8$ only for $l\geq 10$, the new version is able to compute all values for $l\geq 7$.  

Our program is based on the algorithm outlined in \cite{GOconj} and also uses ideas by B\"urgstein and Hesselink \cite{BUHE} as well as Vera-L\'opez and Arregi \cite{VLAR}. 

We now give a brief outline of the structure of the paper. In Section \ref{sec:theor} we give a summary of the theoretical results that the program is based on. The algorithm, with an emphasis on the improvements that have been made, is described in Section \ref{sec:algor}
and the results of our calculations are presented in Section \ref{sec:resul}. Finally in Section \ref{sec:coeff}, we prove explicit formulas for the coefficients of $k(U(q))$ of degrees zero, one and two, assuming that $k(U(q))$ is a polynomial in $q-1$.

\section{Theoretical background}
\label{sec:theor}

Let $G$ be a connected reductive algebraic group, defined and split over the finite field $\F_q$ with $q$ elements, where $q$ is a power of a prime $p$. Assume that $p$ is a good prime for $G$ and let $K$ be the algebraic closure of $\F_q$. We identify $G$ with its group of points over $K$ and write $G(q)$ for the group of $\F_q$-rational points of $G$. Let $B$ be a Borel subgroup of $G$ defined over $\F_q$, containing a maximal torus $T$ defined over $\F_q$. The unipotent radical $U$ of $B$ is also defined over $\F_q$ and the group $U(q)$ of $\F_q$-rational points of $U$ is a Sylow $p$-subgroup of $G(q)$. We write $\u$ for the Lie algebra of $U$ and $\u(q)$ for its space of $\mathbb F_q$-rational points.

Below we recall some of the results from \cite{GOconj} and \cite{GOzeta} on which our algorithm for calculating $k(U(q))$ is based. We note that \cite[Thm.\ 1.1]{GOspring} implies that some of the results in \cite{GOconj} hold in greater generality than stated there and allow ourselves to give the more general statements below. Thanks to \cite[Prop.\ 6.2]{GOconj}, we know that the conjugacy classes of $U(q)$ are in bijective correspondence with the adjoint $U(q)$-orbits in $\u(q)$, so we are henceforth primarily concerned with these orbits.

Let $\Phi$ be the root system determined by $G$ and $T$ and let $\Phi^+$ be the set of positive roots determined by $B$. Denote by $\preccurlyeq$ the partial order on $\Phi$ determined by $\Phi^+$. Let $N$ be the cardinality of $\Phi^+$ and fix an enumeration of $\Phi^+=\{\beta_1,\dots,\beta_N\}$ such that $i\leq j$ whenever $\beta_i\preccurlyeq\beta_j$. Let $\g_{\beta}$ be the root subspace of $\g$ for $\beta\in\Phi$, and fix a Chevalley basis $\{e_{\beta} \mid \beta\in\Phi^+\}$ for $\u$ with $e_{\beta}\in\u(q)$ for each $\beta\in\Phi^+$. For $0\leq i\leq N$, we define
$$ \m_i=\bigoplus^{N}_{j=i+1}\g_{\beta_j}. $$ 

For $x\in\u$, denote by $x+K e_{\beta_i}+\m_i$ the coset $\{x+a_i e_{\beta_i}+\m_i\mid a_i\in K\}$ in $\u/\m_i$. We have the following dichotomy given by \cite[Lem.\ 5.1]{GOconj}: 

\begin{itemize}
\item[(I)] Either all elements of $x+K e_{\beta_i}+\m_i$ are conjugate in $\u/\m_i$ by $U$ (in which case we call $i$ an \textit{inert point} of $x$), or
\item[(R)] no two elements of $x+K e_{\beta_i}+\m_i$ are conjugate in $\u/\m_i$ by $U$ (in which case we call $i$ a \textit{ramification point} of $x$). 
\end{itemize} 

An element $x=\sum_{i=1}^N a_i e_{\beta_i}\in\u$ is said to be the \textit{minimal representative} of its $U$-orbit if $a_i=0$ whenever $i$ is an inert point of $x$. It follows from \cite[Prop.\ 5.4 and Lem.\ 5.5]{GOconj} that each $U$-orbit in $\u$ contains a unique minimal representative. For our algorithm we use the following characterization of minimal representatives, which follows immediately from the discussion above.

\begin{lem} 
\label{lem:repsandcents}
Let $x=\sum^{N}_{j=1}a_j e_{\beta_j}\in\u$.   
Then $x$ is the minimal representative of its $U$-orbit if $a_i=0$ whenever
$\dim\c_\u(x+\m_i)=\dim\c_\u(x+\m_{i-1})-1$. 
\end{lem}

A crucial theoretical result is that one can view the minimal representatives as elements of certain quasi-affine varieties.  Let $c\in\{\I,\R\}^N$ and define 
$$ \u_c=\{x\in\u\mid\text{$i$ is an inert point of $x$ if and only if $c_i=\I$}\} $$ 
and 
$$ X_c=\left\{\sum_{i=1}^N a_i e_{\beta_i}\in\u_c \:\: \vline \:\: \text{$a_i=0$ if $c_i=\I$}\right\}. $$ 
Then, as stated in \cite[Lem.\ 4.2]{GOzeta}, we have that $X_c$ is a locally closed subvariety
of $\u$, and the adjoint $U$-orbits in $\u$ are in bijection with the points
of $\bigcup_{c\in\{\I,\R\}^N}X_c$. Moreover, this implies that the $U(q)$-orbits in $\u(q)$ are parameterized by the $X_c(q)=X_c\cap\u(q)$.

In the next section we describe our algorithm for calculating $k(U(q))$. The idea of the
algorithm is to determine a decomposition of the varieties $X_c$ as a disjoint union 
of locally closed subvarieties, where these subvarieties are given by the vanishing and non-vanishing of explicitly determined polynomials. The key step in calculating these polynomials is to use the characterization of minimal representatives given by Lemma \ref{lem:repsandcents}. The algorithm then proceeds to determine the number of $\F_q$-rational points in these subvarieties.

\section{The algorithm}
\label{sec:algor}

In this section we describe our algorithm for calculating a parameterization of the adjoint $U$-orbits in $\u$ and then determining $k(U(q))$. In particular, we outline the improvements to the algorithm from \cite{GOROconj} that have enabled us to calculate $k(U(q))$ for Chevalley groups of rank six and seven. The proof given in \cite[\S3]{GOROconj} that the algorithm correctly determines all minimal representatives of $U(q)$-orbits in $\u(q)$ remains valid here, so we only concentrate on explaining the algorithm. 

Before going into some details, we give an overview of how the algorithm works. As mentioned in the previous section, it is based on explicitly determining the varieties $X_c$ for $c\in\{\I,\R\}^N$. 

In fact, first we want to generalize our notation. Let $c\in\{\I,\R_0,\R_\nn\}^i$ for some $i=1,\dots,N$. Then we define
$$ \u_c=\{x+\m_i\in\u_i \mid \text{$j$ is an inert point of $x$ if and only if $c_j=\I$}\} $$
and
$$ X_c=\left\{\sum_{j=1}^i a_j e_{\beta_j}+\m_i\in\u_c \:\: \vline \:\: \text{$a_j=0$ if and only if $c_j\in\{\I,\R_0\}$}\right\}. $$
We define $m_c$ to be the number of $j$ with $c_j=\R_n$. Denote by $1\leq k_1<\dots<k_{m_c}\leq N$ the indices with $c_{k_j}=\R_n$ and define $\beta_{c,j}=\beta_{k_j}$. Then we can write each element of $X_c$ in the form $x_c(a)+\m_i:=\sum^{m_c}_{j=1}a_j e_{\beta_{c,j}}+\m_i$, where $a\in K^{m_c}$. Thus we can canonically identify $X_c$ with a subvariety of $(K^\times)^{m_c}$.

Each $X_c$ can be written as a disjoint union $X_c=\bigcup_{\iota=1}^{l_c} X_c^\iota$, where
$$ X_c^\iota:=\{x_c(a)+\m_i\in X_c \mid \text{$f(a)=0$ for all $f\in A_c^\iota$ and $g(a)\neq 0$ for all $g\in B_c^\iota$}\}, $$
and $A_c^\iota,B_c^\iota\subseteq K[t_1,\dots,t_{m_c}]$. We call the $X_c^\iota$ \textit{families of minimal representatives}. Given $A,B\subseteq K[t_1,\dots,t_{m_c}]$, we define
$$ X_{c,A,B}:=\{x_c(a)+\m_i \in X_c \mid \text{$f(a)=0$ for all $f\in A$ and $g(a)\neq0$ for all $g\in B$}\}. $$
Thus we have $X_c=\bigcup_{\iota=1}^{l_c} X_{c,A_c^\iota,B_c^\iota}$. Note that $X_c$ might be a single family (and often is).

The goal of the algorithm is to calculate $A_c^\iota,B_c^\iota\subseteq K[t_1,\dots,t_{m_c}]$ as above for each $c\in\{\I,\R_0,\R_\nn\}^N$.  Then the $U$-orbits in $\u$ are parameterized by the union of all the $X_c^\iota$ and $k(U(q))$ is given by the sum of the $|X_c^\iota(q)|$.

During the main loop of the algorithm we are considering $c \in \{\I,\R_0,\R_\nn\}^{i-1}$ and looking at a family $X_{c,A,B}$. We also have a ``stack'' of other families that we will consider later: the algorithm is a \textit{depth-first backtrack algorithm} which calculates each family to the end, before considering the next family from the stack.

The key step involves determining the variety $X_{c,A,B}^{+1}$, which consists of elements of the form $x_c(a)+b e_{\beta_i}+\m_i$ such that $x_c(a)+\m_{i-1}\in X_{c,A,B}$ and $x_c(a)+b e_{\beta_i}+\m_i$ is the minimal representative of its $U$-orbit. To do this we first calculate $\dim\c_\u(x_c(a)+\m_i)$ so we can apply Lemma \ref{lem:repsandcents} to determine conditions on $a$ for whether $i$ is an inert or ramification point of $x_c(a)$. Then the main objective is to decompose $X_{c,A,B}^{+1}$ into families. The algorithm ends up with triples $(c^\iota,A^\iota,B^\iota)$ for $\iota=1,\dots,k_c$, $c^\iota=(c,c_i)$ with $c_i\in\{\I,\R_0,\R_\nn\}$ and $A^\iota,B^\iota\subseteq\Z[t_1,\dots,t_{m_{c^\iota}}]$ with $A\subseteq A^\iota$ and $B\subseteq B^\iota$. Moreover, these triples are such that we have a disjoint union
\begin{equation} 
\label{e:split}
X_{c,A,B}^{+1}=\bigcup_{\iota=1}^{k_c} X_{c^\iota,A^\iota,B^\iota}.
\end{equation}
The determination of these families is sometimes easy, for example often we have that $i$ is an inert point for all $x_c(a)$ with $x_c(a)+\m_{i-1}\in X_{c,A,B}$. However, determining $(c^\iota,A^\iota,B^\iota)$ may require some complicated analysis of polynomials. A major improvement in the level of analysis applied here is one of the significant additions to the previous algorithm from \cite{GOROconj}.  

The algorithm then proceeds by considering $(c^1,A^1,B^1)$ and adding $(c^\iota,A^\iota,B^\iota)$ for $\iota = 2,\dots,k_c$ to the stack. 

After all of the families $X_{c,A,B}$ have been processed, the number of $\F_q$-rational points $|X_{c,A,B}(q)|$ is determined. These numbers are summed together to calculate $k(U(q))$. Considerable improvement to the processes for calculating $|X_{c,A,B}(q)|$ are made in this new algorithm. In particular, the algorithm aims to ensure that the polynomials in the sets $A$ and $B$ are linear in one indeterminate, which enables the calculation to be made.

We proceed to give a more detailed description of the algorithm.  First we need to give
some notation that is required in this explanation.

Let $\g$ be the Lie algebra of $G$ and $\g_\C$ be the Lie algebra over $\C$ of the same type. Fix a Chevalley basis of $\g_\C$ and denote by $\g_\Z$ the $\Z$-lattice spanned by this Chevalley basis, so that $\g\cong K\otimes_\Z\g_\Z$. Define
$$ \tilde \u:= \Z[t_1,\dots,t_m] \otimes_\Z  \u_\Z, $$
where $m:=\maxim\{m_c\mid c\in\{\I,\R_0,\R_n\}^N,~X_c\neq\varnothing\}$. We allow ourselves to view $e_{\beta_j}\in\tilde\u$, where $e_{\beta_1},\dots,e_{\beta_N}$ are elements of the Chevalley basis of in $\u_\Z$. For $i\leq N$ and $c\in\{\I,\R_0,\R_n\}^i$, we define
$$ x_c(t):=\sum^{m_c}_{j=1}t_j e_{\beta_{c,j}}\in\tilde\u. $$
For $a=(a_1,\dots,a_{m_c})\in(K^\times)^{m_c}$, we define $x_c(a)\in\u$ by substituting $t_j =a_j$ in $x_c(t)$. Let $y_1,\dots,y_N$ be indeterminates. It is shown in \cite[\S3]{GOROconj} that $\sum_{j=1}^N y_j e_{\beta_j}\in\c_\u(x_c(a))$ if and only if $(y_1,\dots,y_N)$ is a solution of a certain system of linear equations
$$ \sum^{N}_{k=1}P_{jk}^c(a)y_k=0, $$ 
where $P_{jk}^c(t)\in\Z[t]$ are linear polynomials determined by the Chevalley commutator relations. To find out for which $a$ we have $\dim\c_\u(x_c(a)+\m_i)<\dim \c_\u(x_c(a)+\m_{i-1})$, so that Lemma \ref{lem:repsandcents} can be applied, one has to check for which $a$ the rank of the matrix $(P_{jk}^c(a))_{j,k}\in\Mat_{(i-1)\times N}(\Z)$
increases when one appends the $i$-th row. 

Now we give the data that the algorithm is holding at any point during a run. Note that the first three data elements $c$, $A$ and $B$ uniquely determine a family $X_c^\iota$ which we also denote by $X_{c,A,B}$ as above. We give some explanation of the meaning of the data here, but parts can only be fully understood once we have described the algorithm. We use speech marks to identify terminology that has not been explained. 

\begin{itemize}
\item The tuple $c\in\{\I,\R_0,\R_n\}^i$ which determines $x_c(t)\in\tilde\u_i$.
\item The set $A$ of polynomials in $\Z[t]$ which vanish on $X_c^\iota$.
\item The tuple $B$ of polynomials in $\Z[t]$ which have no roots in $X_c^\iota$.
\item For each $f\in A\cup B$, we have associated $\sigma(f)$, which is either equal to $0$
or an indeterminate $t_j$ in which $f$ is linear.
\item The matrix $Q(t)\in\Mat_{i\times N}(\Z[t])$ which comprises the first $i$ rows of $(P_{jk}^c(t))_{j,k}$ in row-reduced form.
\item The tuple $\pi$ containing the ``pivots'' used for the first $i$ row-reductions of $Q(t)$.
\item The stack $S:=\{(c,A,B,\pi,\sigma,Q(t))\}$, an ordered subset of
$$ \bigcup_{i=1}^N\{\I,\R_0,\R_n\}^i\times\mathcal{P}(\Z[t])\times\bigcup_{i\in\N}\Z[t]^i
\times\bigcup_{i=1}^N\{0,1,\dots,N\}^i\times\bigcup_{i\in\N}\{t_1,\dots,t_m\}^i\times
\Mat_{i \times N}(\Z[t]) $$ 
containing information about families that the program has not processed yet.
\item The ``good-families'' set $\gamma$, which contains for each family already processed
enough data from which one can recover the number of $\F_q$-rational points in the family.
\item The ``bad output'' $O:=\{(c,A,B)\}$, a subset of
$$ \{\I,\R_0,\R_n\}^N\times\mathcal{P}(\Z[t])\times\bigcup_{i\in\N}\Z[t]^i $$ 
containing sufficient information about each ``bad'' family.
\end{itemize}

At the beginning of the program, we have the configuration

\begin{itemize}
\item $c:=(\R_n)$,
\item $A:=\varnothing$,
\item $B:=\varnothing$,
\item $\sigma:=\varnothing$,
\item $\pi:=(0)$,
\item $Q(t):=0\in\Mat_{1\times N}(\Z[t])$,
\item $O:=\varnothing$, 
\item $\gamma:=\varnothing$, and
\item $S:=\{(\R_0),\varnothing,\varnothing,(0),\varnothing,0\}$. 

\end{itemize}

The main loop in the algorithm is explained as follows. At each point we are considering a family $X_c^\iota=X_{c,A,B}$ as above. In the explanation below we sometimes speak of \textit{relevant $a$}, by which we mean $a$ such that $x_c(a)\in X_{c,A,B}$.

\noindent \textit{Case 1:} If the length $i-1$ of $c$ is smaller than $N$, then we generate the $i$th row of the matrix $(P_{jk}^c(t))_{j,k}$ and append it to $Q(t)$. Then the following operations are applied to row reduce the $i$-th row. For all $1\leq j\leq i-1$ with $\pi_j\neq 0$, we modify the $i$th row $Q_i(t)$ of $Q(t)$ by setting
$$ Q_i(t):=Q_i(t)\frac{Q_{j,\pi_j}(t)}{\gcd(Q_{i,\pi_j}(t),Q_{j,\pi_j}(t))}-
Q_j(t)\frac{Q_{i,\pi_j}(t)}{\gcd(Q_{i,\pi_j}(t),Q_{j,\pi_j}(t))} $$ 
Note that this leads to $Q_{i,\pi_j}(t)=0$ for all $j$ with $\pi_j\neq 0$. Here the {\sf SINGULAR} \cite{SIN} interface for {\sf GAP} is used to calculate the greatest common divisors. 

The next step depends on the set $L_i$ of non-zero polynomials in $Q_i(t)$ which are not divisible by any polynomial in $A$.

\noindent \textit{Case 1a: $L_i=\varnothing$}: In this case, we have that $i$ is a ramification point of $x_c(a)$ for all relevant $a$, and we set
\begin{itemize}
\item $\pi:=(\pi,0)$,
\item $c:=(c,\R_n)$, and
\item $S:=S\cup\{((c,\R_0),A,B,\pi,\sigma,Q(t))\}$.
\end{itemize}

\noindent \textit{Case 1b: There exists $Q_{i,l}(t) \in L_i$ such that $Q_{i,l}$ is a monomial or divides an $f\in B$}: In this case, we have that $i$ is a inert point of $x_c(a)$ for all relevant $a$, and we set \begin{itemize}
\item $\pi:=(\pi,l)$, and
\item $c:=(c,\I)$.
\end{itemize}

\noindent \textit{Case 1c: $L_i\neq\varnothing$, but no $Q_{i,l}(t)$ as in Case 1b exists}: In this case $i$ can be either an inert point or a ramification point of $x_c(a)$ for relevant $a$. Here we pick a $Q_{i,l}(t)\in L_i$ that is minimal with respect to a total order on $\Z[t]$, comparing first the number of terms of two polynomials, then their degrees and finally their leading coefficients. Then we apply some new subroutines, the \textit{polynomial-resolving subroutine} and the \textit{stack-generating subroutine} and update the data as specified by these subroutines.  

These subroutines are a substantial improvement to the algorithm from \cite{GOROconj}. Their aim is to determine the triples $(c^\iota,A^\iota,B^\iota)$ for $\iota = 1,\dots,k_c$ mentioned above such that \eqref{e:split} holds. 

The algorithm aims to construct the sets $A^\iota$ and $B^\iota$ in some way so that each polynomial in $C^\iota:=(A^\iota\setminus A)\cup(B^\iota\setminus B)$ is linear in one of the indeterminates $t_j$.
Often the elements of $C^\iota$ are irreducible factors of $Q_{i,l}(t)$. Though the situation can get considerably more complicated: when a polynomial $f=h_1t_k+h_2$ is linear in the indeterminate $t_k$, where $h_1,h_2\in\Z[t_1,\dots,t_m]$ are polynomials not involving $t_k$, then it is also necessary to consider when the polynomials $h_1$ and $h_2$ give zero values. The {\sf SINGULAR} \cite{SIN} interface for {\sf GAP} is also used in these processes.
    
The variable $\sigma$ is used to record which indeterminate a polynomial in $A$ or $B$ is linear in. So if we have found that $f=h_1 t_k+h_2$ is linear in $t_k$, then we set $\sigma(f)=t_k$. If $f$ is not linear in any indeterminate, then we set $\sigma(f)=0$. If there is more than one such $t_k$, the program choses the $t_k$ which is most ``suitable'' for subsequent calculations. 

Often when we have a polynomial $f\in A$ which is linear in an indeterminate we perform a substitution to reduce the number of indeterminates. This is done when $t_k$ appears linearly in some $f=h_2 t_k-h_1\in A$. Then we substitute $t_k$ for $\frac{h_1}{h_2}$ in $Q(t)$ as well as in all other elements of $A\cup B$.

A ``trick'' that the program sometimes uses in the polynomial resolving subroutine is to make a linear change of coordinates in the indeterminates so that a polynomial becomes linear in an indeterminate. For example, there might be a polynomial of the form $f(t) = (t_1+t_2)^2+t_1 \in A$. Since neither $t_1$ nor $t_2$ appear linearly in $f(t)$, we cannot solve for either of them. By introducing a new variable $z_2:=t_1+t_2$ and replacing $f(t)$ by $z_2^2+t_1$, we are able to solve for $t_1$ and then make a substitution. Implementing routines to look for such substitutions was a huge challenge, and then it also involved updating other parts of data accordingly.

When these two subroutines are complete the algorithm has calculated the 4-tuples \linebreak $(c^\iota,A^\iota,B^\iota,\sigma^\iota)$. Then the data is updated as follows.
\begin{itemize}
\item $\pi:=(\pi,l)$,
\item $c:=c^1=(c,\I)$,
\item $A:=A^1=A$,
\item $B:=B^1$,
\item $\sigma:=\sigma^1$, and
\item $S:=S\cup\{(c^\iota,A^\iota,B^\iota,\pi,\sigma^\iota,Q(t)) \mid \iota=2,\dots,k_c\}$.
\end{itemize}

\noindent \textit{Case 2:} If $c$ has length $N$, then we have determined the family $X_{c,A,B}$. The program now applies a subroutine on $(c,A,B,\sigma)$ to attempt to calculate the number of $\F_q$-rational points of the family $X_{c,A,B}$. This subroutine is called the \textit{nice conditions subroutine}; it constitutes a major improvement to the algorithm from \cite{GOROconj}.

The first step in the nice conditions subroutine involves checking whether each polynomial $f\in A\cup B$
is linear in some $t_j$, i.e.\ $f=h_1 t_j+h_2$, where $h_1,h_2\in K[t_1,\dots,t_{m_c}]$ not involving $t_j$. This is first done using $\sigma$, but further checks are made to see if elements in $A\cup B$ have become linear in an indeterminate as a consequence of substitutions being made at some point during the run. Then the values $a$ for which $f(a)=0$ are given by $a_j=-\frac{h_2(a)}{h_1(a)}$, so
it is relatively straightforward to count the number of $a$ for which $f(a)$ is zero or non-zero. If this can be applied to all the $f\in A\cup B$, then $|X_c(q)|$ can be calculated. However, a great deal of care
needs to be taken here as there are a variety of potential complications, for example $h_1(a)$ or $h_2(a)$ could be zero. Also there will be a number of dependencies between the conditions calculated for different $f\in A\cup B$, so the algorithm is required to make many checks before it is able to complete the calculation. If the subroutine is successful, then it outputs $d=m_c-|A|$ and a tuple $n(c,\sigma)$ which contains information about the number of indeterminates that can take $q-1-j$ values, for each possible $j$.  The variables $d$ and $n(c,\sigma)$ are later used to calculate $|X_{c,A,B}(q)|$ as explained below. 

The algorithm updates the data depending on the outcome of the nice conditions subroutine.

\begin{enumerate}
\item If the subroutine fails to calculate the number of $\F_q$-rational points in $X_{c,A,B}$, then it sets $O:=O\cup\{(c,A,B)\}$. We name such families \textit{bad families}, as in this case this output would need to be analyzed by hand to determine $|X_{c,A,B}(q)|$. For the cases, where we have run the algorithm, we end up with $O=\varnothing$, which is a significant betterment on the algorithm from \cite{GOROconj}.
\item If the subroutine was successful in calculating $|X_{c,A,B}(q)|$, then we call the family a \textit{good family.} In this case $\gamma$ is updated by setting $\gamma:=\gamma\cup\{(d,n(c,\sigma))\}$,
where $d=m_c-|A|$ and $n(c,\sigma)$ are as above. 
\end{enumerate}

The algorithm proceeds by updating the data from the stack as follows. If $S=\varnothing$, then the main loop terminates. Else we update the variables:
\begin{itemize}
\item $(c,A,B,\pi,\sigma,Q(t))=:\Top(S)$,
\item $S:=S\setminus\{\Top(S)\}$,
\end{itemize}
where $\Top(S)$ is the last element in the ordered set $S$.

After the end of the main loop, the data from $\gamma$ is used to calculate a polynomial $Z(q)$. This polynomial is the sum of $|X_{c,A,B}(q)|$, where $(c,A,B)$ runs over all families of minimal representatives for which the program calculated $|X_{c,A,B}(q)|$. When $(c,A,B)$ corresponds to the element $(d,n(c,\sigma))\in\gamma$, we have
$$ |X_{c,A,B}(q)|=(q-1)^{d-\sum^{|\sigma|}_{j=1}n(c,\sigma)_j}\prod^{|\sigma|}_{j=1}(q-1-j)^{n(c,\sigma)_j}. $$
We note that it turns out that in most cases $A=B=\varnothing$, so the cardinality is a given by $|X_{c,A,B}(q)|=(q-1)^{m_c}$; more complicated situations occur rarely when $G$ has small rank, but with increasing frequency for higher ranks. If $O=\varnothing$, then $Z(q)=k(U(q))$. Otherwise one would have to calculate the number of the $\F_q$-rational points of the families in $O$ by hand. As already mentioned, such hand calculations are not required for the cases on which we have run the algorithm.

We remark that the usage of $\g_\Z$ makes it necessary to be careful about implicit divisions occurring during the calculations. Therefore, the program records the primes that occur in the numerator or denominator of any coefficient of a polynomial from $A$ or $B$. These primes, for which our results may not be valid, are then returned at the end. Fortunately, for all calculations that we have completed, these primes were bad primes.

Finally we mention that the algorithm attempts to normalize coefficients when possible to reduce the number of indeterminates required. The maximal torus $T$ acts on the sets of minimal representatives and can be used to normalize certain coefficients of $x_c(t)$ to be equal to 1 as explained in \cite[\S3]{GOROconj}. Let $c\in\{\I,\R_0,\R_n\}^N$ and let $J$ be a linearly independent subset of $\{\beta_{c,1},\dots, \beta_{c,m_c}\}$. Then we can find for any $a\in(K^\times)^{m_c}$ an element $b\in(K^\times)^{m_c}$ such that $b_j=1$ if $\beta_{c,j}\in J$ and $x_c(a)=t\cdot x_c(b)$ for some $t\in T$. This implies that for $i\in\{1,\dots,N\}$, we have that $i$ is an inert point of $x_c(a)$ if and only if it is an inert point of $x_c(b)$. This trick is useful since it reduces the amount of indeterminates which arise in the program, though we have to take care here: it may be the case that $x_c(a)$ and $x_c(b)$ are not conjugate by
an element of $T(q)$. The algorithm has a routine to normalize coefficients to be equal to $1$ as above,
as long as this is possible by elements in $U(q)$.

We also remark here that for $G(q)$ with root system of type $B_7$ and $C_7$, it turned out that there are situations where, surprisingly, the normalization of certain coefficients leads to more complicated polynomials in the sets $A$ and $B$. This meant that we sometimes had to manually override some normalizations.

\section{The results}
\label{sec:resul}

Table \ref{tab:results} contains the values of $k(U(q))$ for $G(q)$ simple of rank at most 7, except $E_7$, written as polynomials in $v:=q-1$. The polynomials up to rank 5 were already calculated in \cite{GOROconj}, while the polynomials for $G(q)$ of type $A_r$, $r\leq 12$, were given in \cite{VLAR}. The newly obtained polynomials are colored red within the tables.

\begin{table}[!b]
\renewcommand{\arraystretch}{1.4}
\begin{tabular}{|l|l|l|}
\hline
$G(q)$    & $k(U(q))$ \\ 
\hline
$A_1$ & 
$v+1$  \\ 
\hline
$A_2$ & 
$v^2+3v+1$  \\
$B_2$ & 
$2v^2+4v+1$ \\ 
$G_2$ & 
$v^3+5v^2+6v+1$ \\
\hline
$A_3$ & 
$2v^3+7v^2+6v+1$ \\
$B_3,C_3$ & 
$v^4+8v^3+16v^2+9v+1$ \\
\hline 
$A_4$ & 
$5v^4+20v^3+25v^2+10v+1$ \\
$B_4,C_4$ & 
$v^6+11v^5+48v^4+88v^3+64v^2+16v+1$ \\
$D_4$ & 
$2v^5+15v^4+36v^3+34v^2+12v+1$ \\
$F_4$ & 
$v^8+9v^7+40v^6+124v^5+256v^4+288v^3+140v^2+24v+1$ \\
\hline
$A_5$ & 
$v^6+18v^5+70v^4+105v^3+64v^2+15v+1$ \\
$B_5,C_5$ & 
$2v^8+24v^7+132v^6+395v^5+630v^4+500v^3+180v^2+25v+1$ \\
$D_5$ & 
$2v^7+22v^6+106v^5+235v^4+240v^3+110v^2+20v+1$ \\
\hline
$A_6$ & 
$8v^7+84v^6+301v^5+490v^4+385v^3+140v^2+21v+1$ \\ 
\textcolor{red}{$B_6,C_6$} & 
\textcolor{red}{$v^{11}+15v^{10}+112v^9+547v^8+1845v^7+4121v^6+5701v^5+4560v^4+1960v^3$} \\
~ & 
\textcolor{red}{$+~410v^2+36v+1$} \\
\textcolor{red}{$D_6$} & 
\textcolor{red}{$v^{10}+13v^9+87v^8+393v^7+1157v^6+2032v^5+2005v^4+1060v^3+275v^2$} \\
~ & 
\textcolor{red}{$+~30v+1$} \\
\textcolor{red}{$E_6$} & 
\textcolor{red}{$v^{11}+12v^{10}+75v^9+353v^8+1286v^7+3178v^6+4770v^5+4035v^4+1800v^3$} \\
~ & 
\textcolor{red}{$+~390v^2+36v+1$} \\
\hline
$A_7$ & 
$4v^9+74v^8+496v^7+1568v^6+2604v^5+2345v^4+1120v^3+266v^2+28v+1$ \\
\textcolor{red}{$B_7,C_7$} & 
\textcolor{red}{$v^{14}+18v^{13}+158v^{12}+899v^{11}+3740v^{10}+11985v^9+29328v^8+52055v^7$} \\
~ & 
\textcolor{red}{$+~62930v^6+48797v^5+22855v^4+6020v^3+812v^2+49v+1$} \\
\textcolor{red}{$D_7$} & 
\textcolor{red}{$4v^{12}+59v^{11}+417v^{10}+1913v^9+6256v^8+14289v^7+21497v^6+20188v^5$} \\
~ & 
\textcolor{red}{$+~11305v^4+3570v^3+581v^2+42v+1$} \\
\hline
\end{tabular}
\vspace{1mm}
\caption{$k(U(q))$ as polynomials in $v=q-1$.}
\label{tab:results}
\end{table}

The most noteworthy observation is that the polynomials for $B_r$ and $C_r$ coincide for fixed $r$. This has already been noticed for $r\leq 5$, and the equality still holds for $r=6,7$. 

Another phenomenon that was already observed in \cite{VLAR} for type $A_r$ and in \cite{GOROconj} for rank at most 5 is that the coefficients of $k(U(q))$ written as polynomials in $v$ are non-negative integers. A heuristic idea why this may be the case was given in \textit{loc.\ cit.}  

It is noteworthy that the constant coefficient of all calculated polynomials is equal to one. We explain this and prove explicit formulas for the coefficients of degrees one and two in the next section. 

Eamonn O'Brien used the $p$-group conjugacy algorithms available in {\sf MAGMA} \cite{MAG} to confirm the values for $k(U(p))$ in each of the cases in Table 1 for $p$ the smallest good prime for $G$.

We also used a modification of our program to calculate the number $k(U(q),U^{(l)}(q))$ of $U(q)$-conjugacy classes in the $l$-th term of the descending central series of $U(q)$, for certain groups and $l\in\N$. Here we also see that these numbers are given by polynomials in $v$ with non-negative integer coefficients. The results of these calculations are presented in Table \ref{tab:lowcs}.    

\begin{table}[!h]
\renewcommand{\arraystretch}{1.4}
\begin{tabular}{|l|l|l|}
\hline
$G(q)$ & $l$ & $k(U(q),U^{(l)}(q))$ \\ 
\hline
$E_7$ & \textcolor{red}{2} & \textcolor{red}{$v^{14}+14v^{13}+92v^{12}+380v^{11}+1128v^{10}+2675v^9+5694v^8+11565v^7$} \\
~ & ~ & \textcolor{red}{$+19486v^6+21745v^5+13976v^4+4724v^3+755v^2+50v+1$} \\
~ & \textcolor{red}{3} & \textcolor{red}{$3v^{10}+37v^9+253v^8+1193v^7+3767v^6+6724v^5+6194v^4+2798v^3+560v^2$} \\
~ & ~ & \textcolor{red}{$+44v+1$} \\
~ & 4 & $v^9+13v^8+94v^7+512v^6+1600v^5+2312v^4+1499v^3+395v^2+38v+1$ \\
\hline
$E_8$ & \textcolor{red}{7} & \textcolor{red}{$2v^{13}+28v^{12}+188v^{11}+822v^{10}+2838v^9+8987v^8+25419v^7+51513v^6$} \\
~ & ~ & \textcolor{red}{$+60889v^5+37867v^4+11140v^3+1428v^2+70v+1$} \\
~ & \textcolor{red}{8} & \textcolor{red}{$v^{12}+14v^{11}+94v^{10}+449v^9+1830v^8+6381v^7+16610v^6+25867v^5$} \\
 ~ & ~ & \textcolor{red}{$+20935v^4+7620v^3+1155v^2+64v+1$} \\
~ & \textcolor{red}{9} & \textcolor{red}{$v^{10}+21v^9+199v^8+1125v^7+4228v^6+9382v^5+10568v^4+4955v^3+912v^2$} \\
~ & ~ & \textcolor{red}{$+58v+1$} \\
~ & 10 & $v^9+17v^8+135v^7+719v^6+2568v^5+4652v^4+3014v^3+699v^2+52v+1$ \\
\hline
\end{tabular}
\vspace{1mm}
\caption{$k(U(q),U^{(l)}(q))$ for $E_7$ and $E_8$ as polynomials in $v=q-1$.}
\label{tab:lowcs}
\end{table}

\section{The coefficients of $k(U(q))$ of small degree}
\label{sec:coeff}
 
Assuming that $k(U(q))$ is a polynomial in $v$, we now prove that the coefficients of $k(U(q))$ of degrees zero, one and two can be easily determined based on properties of the root system. We start with an elementary lemma; since we were unable to find a proof in the literature, we give a complete argument.

\begin{lem}
\label{lem:deproots}
Let $\Phi$ be an irreducible root system of rank $r\geq 3$. Let $\alpha$, $\beta$ and $\gamma$ be three pairwise distinct linearly dependent positive roots with $\height\gamma\geq\maxim\{\height\alpha,\height\beta\}$. Then at least one of the following statements is true: 

\begin{itemize}
\item $\beta-\alpha\in\Phi$, 
\item $\gamma-\beta\in\Phi$ and $\gamma-\alpha\in\Phi$, 
\item $\gamma-\alpha\in\Phi$, but $\beta+\gamma-\alpha\notin\Phi$, or  
\item $\gamma-\beta\in\Phi$, but $\alpha+\gamma-\beta\notin\Phi$.
\end{itemize}

\proof Choose an embedding of $\Phi$ into the real vector space $V$. It is a well-known fact (e.g. \cite[Ex.\ 9.7]{HU}) that, if $V'\subseteq V$ is the subspace spanned by $\alpha$ and $\beta$, then $\Phi\cap V'=:\Phi'$ is a root system of rank two. Since $r\geq 3$, we know that $\Phi'$ is a proper root subsystem of $\Phi$. Because two distinct positive roots are linearly independent, $\gamma$ can be written as a linear combination of $\alpha$ and $\beta$, and thus $\gamma$ lies in $\Phi'$. Let $\Phi'^+:=\Phi'\cap\Phi^+$, then the three roots are also positive in $\Phi'$, and $\height\gamma\geq\maxim\{\height\alpha,\height\beta\}$ stays true in $\Phi'$. 

There are three types of root systems of rank two which might appear as proper root subsystems in an irreducible root system ($G_2$ never does). We denote by $\delta$ and $\epsilon$ the simple roots of $\Phi'$ and proceed by case-by-case analysis. 

\textit{$\Phi'$ of type $A_1\times A_1$:} This is not possible, since $\alpha$, $\beta$ and $\gamma$ are pairwise distinct. 

\textit{$\Phi'$ of type $A_2$:} We have a bijection between $\{\alpha,\beta,\gamma\}$ and $\Phi'^+=\{\delta,\epsilon,\delta+\epsilon\}$. Then $\gamma=\delta+\epsilon$ due to the height condition, and then the second statement is true. 

\textit{$\Phi'$ of type $B_2$:} Here $\Phi'^+=\{\delta,\epsilon,\delta+\epsilon,\delta+2\epsilon\}$. If $\gamma=\delta+\epsilon$, then the second statement is true again. The remaining case is when $\gamma=\delta+2\epsilon$. If either $\alpha$ or $\beta$ is $\delta+\epsilon$, then the other one is a simple root and the first statement holds. So suppose that $\{\alpha,\beta\}=\{\delta,\epsilon\}$. Depending on whether $\alpha$ is equal to $\delta$ or $\epsilon$, the third or the fourth statement holds. \hfill $\square$
\end{lem}

\begin{lem}
\label{lem:orbsize}
Let $I$ be a subset of $\Phi^+$, define the span 
$$ X_I:=\left\{\sum_{\beta\in I}a_{\beta} e_{\beta} \:\: \vline \:\: a_{\beta}\in K^\times\right\}\subseteq\u $$ 
and let $\beta_1,\dots,\beta_k$ be linearly independent roots in $I$. Let $\tilde\beta_j$ be the coordinate vector of $\beta_j$ with respect to the base of $\Phi$ determined by $\Phi^+$, for $1\leq j\leq k$. Let $d_1,\dots,d_k$ be the diagonal entries of the Smith normal form of the matrix $(\tilde\beta_1,\dots,\tilde\beta_k)\in\Mat_{r\times k}(\Z)$. 

Then for each $x\in X_I(q)$, the size of the orbit $T(q)\cdot x$ is divisible by $v^k/d$, where $v:=q-1$ and $d:=\prod_{l=1}^{k}\gcd\{d_l,v\}$. 

\proof Since $\beta_1,\dots,\beta_k$ are linearly independent, they form a basis for the sub-lattice $L:=\Z\beta_1+\dots+\Z\beta_k\subseteq\Z\Phi$. The theory of finitely generated abelian groups makes it possible to find a basis $\chi_1,\dots,\chi_r$ of the character group $X(T)$ such that $d_1\chi_1,\dots,d_k\chi_k$ is a basis for $L$. Write $\beta_j=\sum_{l=1}^{k}c_{lj} d_l\chi_l$. Then the matrix $(c_{lj})_{l,j}\in\Mat_{k\times k}(\Z)$ is invertible and we denote by $(a_{jh})_{j,h}\in\Mat_{k\times k}(\Q)$ its inverse. Let $\psi_1,\dots,\psi_r\in X^\vee(T)$ be dual to $\chi_1,\dots,\chi_r$ with respect to the perfect pairing $\langle\cdot,\cdot\rangle$ on $X(T)\times X^\vee(T)$, i.e.\ $\langle\psi_l,\chi_j\rangle=\delta_{lj}$. 

It is known that $\beta_j\circ\psi_l(b)=b^{\langle\psi_l,\beta_j\rangle}$ and $t\cdot e_{\beta}=\beta(t)e_{\beta}$, for $b\in K^\times$ and $t\in T$. By looking at the coefficient of $t\cdot x$ belonging to $e_\beta$, we get that $t=\prod_{l=1}^{r}\psi_l(b_l)\in C_T(x)$ satisfies 
$$ \prod_{l=1}^{k}b_l^{c_{lj}d_l}=1~\text{for all}~1\leq j\leq k. $$ 
Moreover, it follows from
$$ \prod_{j=1}^{k}\left(\prod_{l=1}^{k}t_l^{c_{lj}d_l}\right)^{a_{jh}}=t_h^{d_h} $$ 
that 
$$ C_T(x)\subseteq S:=\left\{\prod_{l=1}^{r}\psi_l(b_l) \:\: \vline \:\: b_l\in K^\times~\text{for all}~1\leq l\leq r,~b_l^{d_l}=1~\text{for all}~1\leq l\leq k\right\}. $$ 
If we take $x\in X_I(q)$ and consider the action of $T(q)$ on $X_I(q)$, then $C_{T(q)}(x)\subseteq S(q)$. Because of the condition on the $b_j$ to be $d_j$-th roots of unity and to lie in $\F_q$ (i.e. to be $v$-th roots of unity as well), the order of $S(q)$ is $v^{r-k}d$, thus the order of $C_{T(q)}(x)$ divides this number. By the orbit stabilizer theorem, $T(q)\cdot x$ has size $v^r/|C_{T(q)}(x)|$, which is divisible by $v^k/d$. \hfill $\square$
\end{lem}

\begin{lem}
\label{lem:transitact}
Let $\beta_1,\dots,\beta_k$ be linearly independent roots in $\Phi^+$. Then: 
\begin{enumerate} 
\item $T$ acts transitively on $X:=\{\sum_{j=1}^{k}a_j e_{\beta_j} \mid a_j\in K^\times\}$, and
\item $\dim\c_\u(x+\m_i)=\dim\c_\u(y+\m_i)$ for all $x,y\in X$, $1\leq i\leq N$. 
\end{enumerate}

\proof (1): We define $\psi_1,\dots,\psi_r\in X^\vee(T)$ and $d_1,\dots,d_k$ as in Lemma \ref{lem:orbsize}. For $x:=\sum_{j=1}^{k}a_j e_{\beta_j}\in X$ and $t=\prod_{l=1}^{r}\psi_l(b_l)$ it follows that $t\cdot x=\sum_{j=1}^{k}b_i^{d_j}a_j e_{\beta_j}$. By taking $b_j$ to be a $j$-th root of unity of $a_j^{-1}$ for all $j$, we see that $x$ lies in the same $T$-orbit as $\sum_{j=1}^{k}e_{\beta_j}$. By transitivity, all $x\in X$ are $T$-conjugate.  

(2): Let $x$ and $y=t\cdot x$ be in $X$. Then also $x+\m_i=t\cdot(y+\m_i)$ for all $1\leq i\leq N$. We get 
\begin{align*}
tC_U(x+\m_i)t^{-1} & =\{\tilde{u}:=tut^{-1}\in U \mid u\cdot(x+\m_i)=x+\m_i\} \\
 & =\{\tilde{u}\in U \mid \tilde{u}t\cdot(x+\m_i)=t\cdot(x+\m_i)\} \\
 & =\{\tilde{u}\in U \mid \tilde{u}\cdot(y+\m_i)\tilde{u}^{-1}=y+\m_i\} \\
 & =C_U(y+\m_i), 
\end{align*} 
so $\dim C_U(x+\m_i)=\dim C_U(y+\m_i)$. Now (2) follows from \cite[Cor.\ 4.3]{GOconj}. \hfill $\square$
\end{lem}

\begin{thm}
\label{thm:degcoeff}
If $k(U(q))$ is given by a polynomial in $v:=q-1$, then the following statements hold: 
\begin{enumerate}
\item The coefficient of degree zero equals 1. 
\item The coefficient of degree one equals $|\Phi^+|$.
\item The coefficient of degree two equals $|\{(\beta_j,\beta_k)\in\Phi^+\times\Phi^+ \mid j<k, ~\beta_k-\beta_j\notin\Phi\}|$. 
\end{enumerate}

\proof (1):  We want to prove that $k(U(q))-1$ is divisible by $v/d$ for some fixed $d\in\N$. First, we note that there is exactly one family of minimal representatives with all coefficients being zero, namely $\{0\}\subseteq\u$. Now, if $X_c$ is a different set of minimal representatives, then there is at least one non-zero coefficient (i.e.\ $m_c>0$). Lemma \ref{lem:orbsize} with $k=1$ and $\beta_1=\beta_c$ yields that there is a $d_{c,v}\in\N$ such that the size of each $T(q)$-orbit on $X_c(q)$ (and thus the cardinality of $X_c(q)$) is divisible by $v/d_{c,v}$. Take 
$$ d_v:=\lcm\{d_{c,v}\mid X_c~\text{set of minimal representatives with}~m_c>0\}.$$ 
Then $v/d_v$ divides $|X_c(q)|$ for all $X_c$ with $m_c>0$. Thus $v/d_v$ divides $k(U(q))-1$. Due to the definition of the $d_{c,v}$ in Lemma \ref{lem:orbsize}, there is a $d\in\N$ such that $d_v$ divides $d$ for all $v\in\N$. Since $k(U(q))$ is a polynomial in $v$, it follows that $v/d$ divides $k(U(q))-1$ as a polynomial. 

Note that for $k=1$ one can give a simpler argument because $d_{c,v}=1$ for all $c,v$. However, the argument above is needed for (2) and (3).   

(2): Let $X_c$ be a set of minimal representatives with $m_c=1$, i.e.\  $X_c\subseteq\{a_j e_{\beta_j} \mid a_j\in K^\times\}$ for some $1\leq j\leq N$. It follows from Lemma \ref{lem:transitact}(1) and \cite[Lem.\ 7.2]{GOconj} that $|X_c(q)|=v$. The number of such sets is $N=|\Phi^+|$. 

Now we consider sets of minimal representatives with $m_c>1$. Since two distinct positive roots are linearly independent, we can use again Lemma \ref{lem:orbsize} (this time with $k=2$) and argue as in (1) that $k(U(q))-|\Phi^+|v-1$ is divisible by $v^2/d$ for some $d\in\N$. 

(3): Consider a set $X_c\subseteq\{a_j e_{\beta_j}+a_k e_{\beta_k} \mid a_j,a_k\in K^\times\}$ of minimal representatives that has two non-zero coefficients, with $j<k$. Lemma \ref{lem:transitact}(2) implies that whether $k$ is an inert point of $x\in X_c$ only depends on $j$ (and not on $a_j$). Thus we can consider $x:=e_{\beta_j}+e_{\beta_k}$. The fact that $k$ is a ramification point of $x$ is equivalent to there being no positive root $\alpha$ such that $[e_{\beta_j},e_{\alpha}]=c e_{\beta_k}$, with $c\neq 0$ (else $\dim\c_{\u}(x+\m_k)=\dim\c_{\u}(x+\m_{k-1})-1$). By Chevalley's commutator formula, $\alpha=\beta_k-\beta_j$, i.e.\ $\beta_k-\beta_j$ must not be a root. 

It follows that there are as many sets $X_c$ of minimal representatives with $m_c=2$ as there are tuples $(\beta_j,\beta_k)$ of positive roots with $j<k$ and $\beta_k-\beta_j\notin\Phi$. Because of Lemma \ref{lem:transitact}(1) and \cite[Lem.\ 7.2]{GOconj}, it follows that $|X_c(q)|=v^2$ for these $X_c$. 

Now, let $X_c$ be a set of minimal representatives with more than two non-zero coefficients, and let the first three of them (with respect to our ordering) belong to $e_{\beta_j}$, $e_{\beta_k}$ and $e_{\beta_l}$. We first want to prove that $\beta_j$, $\beta_k$ and $\beta_l$ must be linearly independent. We can use again Lemma \ref{lem:transitact}(2) and consider the sum $e_{\beta_j}+e_{\beta_k}+e_{\beta_l}$. From the fact that $j$, $k$ and $l$ are ramification points we can deduce the following facts:

\begin{itemize}
\item \textit{$\beta_k-\beta_j$ is not a root.} This follows similarly as in (2), because $k$ is a ramification point. 
\item \textit{Not both $\beta_l-\beta_j$ and $\beta_l-\beta_k$ are roots.} Otherwise, if $\beta_l-\beta_j=:\beta_m$ and $\beta_l-\beta_k=:\beta_n$, then $l$ being a ramification point implies that there is a dependence between the coefficients of $e_{\beta_m}$ and $e_{\beta_n}$ in $\c_{\u}(e_{\beta_j}+e_{\beta_k}+e_{\beta_l})$. This dependence could have only originated from an earlier inert point $s$, and so $[e_{\beta_m},e_{\beta_k}]$ and $[e_{\beta_n},e_{\beta_j}]$ must be non-zero elements from the root space $e_{\beta_s}$. Using again Chevalley's commutator formula, we get that $\beta_m+\beta_k=\beta_j+\beta_n$. Together with $\beta_m+\beta_j=\beta_k+\beta_n$ this leads to $\beta_j=\beta_k$, a contradiction. 
\item \textit{If $\beta_l-\beta_j=:\beta_m$ is a root, then $\beta_k+\beta_m$ is also a root.} Since $l$ is a ramification point in spite of $[e_{\beta_j},e_{\beta_m}]=c e_{\beta_l}$ for some $c\neq0$, the centralizer $\c_{\u}(e_{\beta_j}+e_{\beta_k}+e_{\beta_l})$ must consist of elements with the coefficient of $e_{\beta_m}$ being zero. This must originate from an earlier inert point $s<l$, which means that $[e_{\beta_m},e_{\beta_k}]=c e_{\beta_s}$ for some $c\neq0$. The Chevalley commutator formula yields that $\beta_m+\beta_k=\beta_s$. 
\item \textit{If $\beta_l-\beta_k$ is a root, then $\beta_j+\beta_l-\beta_k$ is also a root.} This follows analogously to the previous statement. 
\end{itemize}

If our group $G$ has rank one or rank two, the statement of the theorem follows from the respective polynomials which are known. If the rank of $G$ is bigger than two, then the linear independence is a direct result of Lemma \ref{lem:deproots}, using contraposition. 

Now we can again use Lemma \ref{lem:orbsize} (with $k=3$) and argue as in (1) and (2) that 
\[
k(U(q))-|\{(\beta_j,\beta_k)\in\Phi^+\times\Phi^+ \mid j<k,~\beta_k-\beta_j\notin\Phi\}|v^2-|\Phi^+|v-1
\]
is divisible by $v^3/d$ for some $d\in\N$. \hfill $\square$
\end{thm}

We summarize the idea behind the proof of the formulas in Theorem \ref{thm:degcoeff}: Suppose that we have a set $J$ of $k$ positive roots in a root system of rank $r\geq k$, where $k\leq 3$. If there exists a family $X_c$ of minimal representatives such that the nonzero coefficients of the elements of $X_c$ correspond to the roots in $J$, then the roots in $J$ must be linearly independent. This statement is trivial for $k=1$ and $k=2$, but some work is required for $k=3$.

\begin{exmp}
\label{exmp:counterex}
It is not possible to use a similar argument in order to determine a formula for the coefficients of degrees three and higher: Let $\Phi$ be of type $C_5$ with basis $\{\alpha_1,\dots,\alpha_5\}$, where $\alpha_5$ is long, and consider $J:=\{\alpha_2,~\alpha_5,~\alpha_3+\alpha_4+\alpha_5,~2\alpha_3+2\alpha_4+\alpha_5\}$. Then the set 
$$X_c:=\{x=a_1 e_{\beta_1}+\dots+a_n e_{\beta_n}\mid a_i\neq 0~\Leftrightarrow~\beta_i\in J\}$$ 
is a family of minimal representatives that occurs for type $C_5$, but the roots in $J$ are not linearly independent any more. So the aforementioned statement does not hold for $k>3$. 
\end{exmp}

\bigskip

\textbf{Acknowledgments:} The research for this paper was carried out in part while the authors were staying at the Mathematical Research Institute Oberwolfach supported by the ``Research in Pairs'' programme. Part of this paper was written during a stay of the second author at the University of Birmingham supported by a stipend of the Ruth and Gert Massenberg Foundation, and also during a visit of the first author to the Ruhr-University Bochum. Finally, we would like to thank Frank L\"ubeck for helpful discussions and Eamonn O'Brien for confirming some of our results for the smallest good primes.

\end{document}